\documentclass[a4paper,12pt]{amsart}
\usepackage{amsfonts}
\usepackage{amssymb}
\usepackage{ifthen}
\usepackage{graphicx}

\nonstopmode \numberwithin{equation}{section}
\newtheorem{thm}{Theorem}[section]

\newtheorem{cor}{Corollary}[section]
\newtheorem{lem}{Lemma}[section]
\newtheorem{prop}{Proposition}[section]

\newtheorem{claim}{Claim}
\newtheorem{conj}[equation]{Conjecture}

\theoremstyle{definition}
\newtheorem{defn}{Definition}[section]
\newtheorem{examp}{Example}[section]
\newtheorem{prob}[equation]{Problem}
\newtheorem{ques}[equation]{Question}
\newtheorem{rem}{Remark}[section]

\newcounter {own}
\def\theown {\thesection       .\arabic{own}}

\newenvironment{pf}[1][]{%
 \vskip 3mm
 \noindent
 \ifthenelse{\equal{#1}{}}%
  {{\slshape Proof. }}%
  {{\slshape #1.} }%
 }%
{\qed\bigskip}

\newcounter{alphabet}
\newcounter{tmp}

\newcounter{alphabet2}

\makeatletter
\newcommand{\Ref}[1]{\@ifundefined{r@#1}{}{\setcounter{tmp}{\ref{#1}}\Alph{tmp}}}
\makeatother
\newenvironment{Lem}[1][]{\refstepcounter{alphabet}%
\bigskip%
\noindent%
{\bf Lemma \Alph{alphabet}}%
{\bf .} \itshape}{\vskip 8pt}

\newcommand{\IC}{{\mathbb C}}

\newcommand{\ID}{{\mathbb D}}

\newcommand{\real}{{\operatorname{Re}\,}}

\def\be{\begin{equation}}
\def\ee{\end{equation}}

\newcommand{\bee}{\begin{enumerate}}
\newcommand{\eee}{\end{enumerate}}

\newcommand{\blem}{\begin{lem}}
\newcommand{\elem}{\end{lem}}
\newcommand{\bthm}{\begin{thm}}
\newcommand{\ethm}{\end{thm}}
\newcommand{\bcor}{\begin{cor}}
\newcommand{\ecor}{\end{cor}}
\newcommand{\beg}{\begin{examp}}
\newcommand{\eeg}{\end{examp}}
\newcommand{\begs}{\begin{examples}}
\newcommand{\eegs}{\end{examples}}
\newcommand{\bdefe}{\begin{defn}}
\newcommand{\edefe}{\end{defn}}
\newcommand{\bprob}{\begin{prob}}
\newcommand{\eprob}{\end{prob}}
\newcommand{\bques}{\begin{ques}}
\newcommand{\eques}{\end{ques}}
\newcommand{\bei}{\begin{itemize}}
\newcommand{\eei}{\end{itemize}}

\newcommand{\bca}{\begin{case}}
\newcommand{\eca}{\end{case}}
\newcommand{\bcl}{\begin{claim}}
\newcommand{\ecl}{\end{claim}}

\newcommand{\bcon}{\begin{conj}}
\newcommand{\econ}{\end{conj}}
\newcommand{\bcons}{\begin{conjs}}
\newcommand{\econs}{\end{conjs}}
\newcommand{\bprop}{\begin{prop}}
\newcommand{\eprop}{\end{prop}}
\newcommand{\br}{\begin{rem}}
\newcommand{\er}{\end{rem}}
\newcommand{\brs}{\begin{rems}}
\newcommand{\ers}{\end{rems}}
\newcommand{\bo}{\begin{obser}}
\newcommand{\eo}{\end{obser}}
\newcommand{\bos}{\begin{obsers}}
\newcommand{\eos}{\end{obsers}}
\newcommand{\bpf}{\begin{pf}}
\newcommand{\epf}{\end{pf}}
\newcommand{\ba}{\begin{array}}
\newcommand{\ea}{\end{array}}
\newcommand{\beq}{\begin{eqnarray}}
\newcommand{\beqq}{\begin{eqnarray*}}
\newcommand{\eeq}{\end{eqnarray}}
\newcommand{\eeqq}{\end{eqnarray*}}

\newcommand{\ds}{\displaystyle}

\begin{document}

\bibliographystyle{amsplain}

\title[On properties of the solutions to the $\alpha$-harmonic equation]
{On properties of the solutions to the $\alpha$-harmonic equation}
\author{Peijin Li}
\address{Peijin Li,
 Department of Mathematics,
Hunan First Normal University, Changsha, Hunan 410205, People's Republic of China}
\email{wokeyi99$@$163.com}

\author{Antti Rasila$^*$}
\address{Antti Rasila,
Department of
Mathematics and Systems Analysis, Aalto University, P. O. Box 11100,
FI-00076 Aalto, Finland.
}
\email{antti.rasila$@$iki.fi}

\author{Zhi-Gang Wang}
\address{Zhi-Gang Wang,
Department of Mathematics,
Hunan First Normal University, Changsha, Hunan 410205, People's Republic of China}
\email{wangmath$@$163.com}

\date{\today}


\thanks{$^*$Corresponding author.}

%\dedicatory{}
\begin{abstract}
The aim of this paper is to establish properties of the solutions to the $\alpha$-harmonic equations: $\Delta_{\alpha}(f(z))=\partial{z}[(1-{|{z}|}^{2})^{-\alpha} \overline{\partial}{z}f](z)=g(z)$, where $g:\overline{\ID}\rightarrow\mathbb{C}$ is a continuous function and $\overline{\ID}$ denotes the closure of the unit disc $\ID$ in the complex plane $\mathbb{C}$. We obtain Schwarz type and Schwarz-Pick type inequalities for the solutions to the $\alpha$-harmonic equation.
%As an application, we get the Lipschitz continuity with respect to the hyperbolic metric of the solutions to $\alpha$-harmonic equations.
In particular, for $g\equiv 0$, the solutions to the above equation are called $\alpha$-harmonic functions. We determine the necessary and sufficient conditions for an analytic function $\psi$ to have the property that $f\circ\psi$ is $\alpha$-harmonic function for any $\alpha$-harmonic function $f$. Furthermore, we discuss the Bergman-type spaces on $\alpha$-harmonic functions.
\end{abstract}

%\thanks{The research was partly supported by NSFs of China (No. 11571216 and No. 11671127)}

\maketitle % \pagestyle{myheadings}
% \markboth{Peijin Li, Antti Rasila and Zhi-Gang Wang}{Some properties of solutions to $\alpha$-harmonic equations}

\section{Introduction}\label{csw-sec1}

Let $\mathbb{C}$ denote the complex plane. For $a\in$ $\mathbb{C}$, let $\mathbb {D}(a,r)=\{z:|z-a|<r\}$ $(r>0)$  and $\mathbb {D}(0,r)=\mathbb {D}_{r}$, $\mathbb {D}=\mathbb {D}_1$ and $\mathbb{T}=\partial\mathbb{D}$, the boundary of $\mathbb{D}$, and $\overline{\ID}=\ID\cup \mathbb{T}$, the closure of $\mathbb{D}$. Furthermore, we denote by $\mathcal{C}^{m}(\Omega)$ the set of all
complex-valued $m$-times continuously differentiable functions from
$\Omega$ into $\mathbb{C}$, where $\Omega$ stands for a subset of
$\mathbb{C}$ and $m\in\mathbb{N}_0:=\mathbb{N}\cup\{0\}$. In
particular, $\mathcal{C}(\Omega):=\mathcal{C}^{0}(\Omega)$ denotes the
set of all continuous functions in $\Omega$. We use $d(z)$ to denote the Euclidean distance from $z$ to $\mathbb{T}$.

\subsection{Distributions}

Let $\Omega$ be an open set in $\mathbb{R}^n$, and let $f\colon \Omega \to \mathbb{R}$ be an infinitely many times
differentiable function. We may write its partial derivatives in the form
$$\partial^{\alpha}f=\partial_1^{\alpha_1}\cdots\partial_n^{\alpha_n}f,$$
where $\alpha=(\alpha_1,\ldots, \alpha_n)$ is a {\it multi-index} and $\partial_j=\partial/\partial x_j$, $j=1, 2, \ldots, n$.
We denote by $C_0^{\infty}(\Omega)$ the space of functions which are infinitely many times differentiable and have a compact support in $\Omega$.

A distribution $f$ in $\Omega$ is a linear form on $C_0^\infty( \Omega)$ such
that for every compact set $K\subset  \Omega$ there exist constants $C$ and $k$ such that
$$|f(\phi)|\leq C\sum_{|\alpha|\leq k}\sup|\partial^{\alpha}\phi|,$$
where $\phi\in C_0^\infty(K)$.
The set of all distributions in $\Omega$ is denoted by $\mathcal{D}'( \Omega)$ (cf. \cite{LH}).

For locally integrable $u$ and $v$ are on $\Omega$ and
$$-\int_{\Omega}u\partial_j \varphi\, dx=\int_{\Omega}v\varphi \,dx,$$
for all $\varphi\in C_0^{\infty}(\Omega)$, we shall say that $\partial_j u=v$ {\it in the distribution sense} (cf. \cite{GG}).

\subsection{The $\alpha$-harmonic equation}

We denote by $\Delta _{\alpha}$ the weighted Laplace operator corresponding to the so-called standard weight $w_{\alpha}=(1-|z|^{2})^{\alpha}$, that is, \beq\label{eq0.1}
\Delta_{\alpha}(f(z))=\partial{z}[(1-{|{z}|}^{2})^{-\alpha} \overline{\partial}{z}f](z)
\eeq
in $\mathbb{D}$, where $\alpha > -1$ (see \cite[Proposition $1.5$]{oA} for the reason for this constraint),
$$\partial{z}=\frac{1}{2}\Big(\frac{\partial }{\partial{x}}-i\frac{\partial }{\partial{y}}\Big)\;\;\mbox{and}\;\;
\overline{\partial}{z}=\frac{1}{2}\Big(\frac{\partial }{\partial{x}}+i\frac{\partial }{\partial{y}}\Big).$$

The weighted Laplacians of the form \eqref{eq0.1} were first systematically studied by %Paul
 Garabedian in \cite{Ga}.
In \cite{oA}, Olofsson and Wittsten introduced the operator $\Delta _{\alpha}$ and gave a counterpart of the classical Poisson integral formula for it.

Let $g\in\mathcal{C}(\mathbb{D})$
and $f\in\mathcal{C}^{2}(\mathbb{D})$. Of particular interest to us is the following {\it inhomogeneous $\alpha$-harmonic equation} in $\mathbb{D}$:
\beq\label{eq1.1}
\Delta_{\alpha}(f)=g.
\eeq
We also consider the associated \emph{Dirichlet boundary value problem} of functions $f$, satisfying the equation \eqref{eq1.1},
\beq\label{eq1.2}
\left\{
\begin{aligned}
\Delta_{\alpha}(f)&=g\;\;\;\;\text{in}\;\mathbb{D},\\
                 f&=f^{\ast}\;\;\text{on}\;\mathbb{T}.
\end{aligned}
\right.
\eeq
Here the boundary data $f^{\ast}$ is a distribution on $\mathbb{T}$, i.e. $f^{\ast}\in \mathcal{D}'(\mathbb{T})$,  and the boundary condition in \eqref{eq1.2} is  understood as $f_{r}\to f^{\ast}\in \mathcal{D}'(\mathbb{T})$ as $r \to 1^{-}$, where

$$
f_{r}(e^{i\theta})=f(re^{i\theta})
$$
for $e^{i\theta}\in\mathbb{T}$ and $r\in[0,1)$.

If $g\equiv0$ in \eqref{eq1.1}, the solutions to \eqref{eq1.1} are said to be {\it$\alpha$-harmonic functions}. Obviously, $\alpha$-harmonicity coincides with harmonicity when $\alpha = 0$. See \cite{DuR} and the references therein for the properties of harmonic mappings.

In \cite{oA}, Olofsson and Wittsten showed that if an $\alpha$-harmonic function $f$ satisfies
$$\lim_{r\to 1^{-}}f_{r}=f^{\ast} \in \mathcal{D}'(\mathbb{T}) \;\; (\alpha > -1),$$ then it has the form of a \emph{Poisson type integral}
\beq\label{eq1.4}
f(z)=\mathcal{P}_{\alpha}[f^{\ast}](z)=\frac{1}{2\pi}\int^{2\pi}_{0}P_{\alpha}(ze^{-i\theta})f^{\ast}(e^{i\theta})d\theta
\eeq in $\mathbb{D}$,
where
$$P_{\alpha}(z)=\frac{(1-|z|^{2})^{\alpha+1}}{(1-z)(1-\overline{z})^{\alpha+1}}$$
is the {\it $\alpha$-harmonic Poisson kernel} in $\ID$. See \cite{LWX} and \cite{olo} for related discussions.

In \cite{B}, Behm found the Green function for $-\overline{\Delta _{\alpha}}$. The {\it $\alpha$-harmonic Green function} $G_{\alpha}$ is given in $\ID$
by
\beq\label{eqGr}
G_{\alpha}(z,w)=-(1-z\overline{w})^{\alpha}h\circ\varphi(z,w),\;\;z\neq w,\eeq
where
\beq\label{eqhs}
h(s)=\int_{0}^{s}\frac{t^{\alpha}}{1-t}dt=\sum^{\infty}_{n=0}\frac{s^{\alpha+1+n}}{\alpha+1+n},\;\;0\leq s<1,\eeq
and
$$\varphi(z,w)=1-\Big|\frac{z-w}{1-\overline{z}w}\Big|^2=\frac{(1-|z|^2)(1-|w|^2)}{|1-z\overline{w}|^2}.$$

For convenience, we let
$$\mathcal{G}[g](z)=\int_{\ID}G_{\alpha}(z, w)g(w)\,d A(w)
$$
and
$$\mathcal{P}[f^{\ast}](z)=\frac{1}{2\pi}\int^{2\pi}_{0}P(ze^{-i\theta})f^{\ast}(e^{i\theta})d\theta,$$
where $dA(w) =(1/\pi)\, dx\,dy$ denotes the normalized area measure in $\ID$ and
$$P(z)=\frac{1-|z|^2}{|1-z|^2}$$
is the {\it Poisson kernel} in $\ID$.

By \cite[Theorem 5.3]{oA} and \cite[Theorem 2]{B}, we see that all
solutions to the $\alpha$-harmonic equation (\ref{eq1.2}) are given by
\beq%\label{eq1.4}
f(z)=\mathcal{P}_{\alpha}[f^{\ast}](z)+\mathcal{G}[g](z).
\eeq

%In the following, we always assume that every $\alpha$-harmonic function has such a representation which plays a key role in the discussions of this paper.

\section{Main results}

\subsection{A Schwarz type lemma}

The classical Schwarz lemma states that an analytic function
$f$ from $\mathbb{D}$ into itself, with $f(0)=0$, satisfies
$|f(z)|\leq|z|$ for all $z\in\mathbb{D}$. This result is a crucial theme in many branches of
 research for more than a hundred years.

Heinz \cite{He} proved the following Schwarz lemma of complex-valued harmonic
functions: If $f$ is a complex-valued harmonic function from
$\mathbb{D}$ into itself with $f(0)=0$, then,
for $z\in\mathbb{D}$, \be\label{eqh}|f(z)|\leq\frac{4}{\pi}\arctan
|z|.\ee

The first purpose of this paper is to consider the results of the above
type for the solutions to the $\alpha$-harmonic equation. Our result is the following:

\begin{thm}\label{thm-1}
Suppose that $g\in\mathcal{C}(\overline{\mathbb{D}})$ and $f^{\ast}\in \mathcal{C}^{1}(\mathbb{T})$. If $f\in\mathcal{C}^{2}(\mathbb{D})$ satisfying \eqref{eq1.2} with $\alpha\geq0$ and $\mathcal{P}_{\alpha}[f^{\ast}](0)=0$, then for $z\in\overline{\ID}$,
\be\label{eqs}
|f(z)|\leq 2^{\alpha}\left[\frac{4}{\pi}\|f^{\ast}\|_{\infty}\arctan|z|+\|g\|_{\infty}(1-|z|^2)^{\alpha+1}\right],
\ee
where
$\|f^{\ast}\|_{\infty}=\sup_{z\in\mathbb{T}}\{|f^{\ast}(z)|\}$, and $\|g\|_{\infty}=\sup_{z\in\mathbb{D}}\{|g(z)|\}$.

Moreover, if we take $\alpha=0$, $g(z)\equiv -C$, where $C$ is a positive constant, and
$$f(z)=C(1-|z|^{2}),$$ then we see that the inequality \eqref{eqs} is sharp.
\end{thm}

Clearly, if $\alpha=0$, $g\equiv 0$ and $f$ maps $\ID$ into itself, then \eqref{eqs} coincides with \eqref{eqh}.

\subsection{A Schwarz-Pick type lemma}
Suppose $f=u+iv$ is in $C^1(\Omega)$, where $\Omega$ is a domain in $\mathbb{C}$ and $u,v$ are real functions. The Jacobian matrix of $f$ at $z$ is denoted by
$$
D_f(z) =
\left(
  \begin{array}{cc}
    u_x & u_y \\
    v_x & v_y \\
  \end{array}
\right).$$
Then
\be\label{qh-1}
\|D_f(z)\|=\sup\{|D_f(z)\varsigma|:\; |\varsigma|=1\}=|f_z(z)|+|f_{\overline z}(z)|
\ee
and
\be\label{qh-2}
l(D_f(z))=\inf\{|D_f(z)\varsigma|:\; |\varsigma|=1\}=\big ||f_z(z)|-|f_{\overline z}(z)|\big|.
\ee

Colonna \cite{Co} obtained a sharp Schwarz-Pick type lemma for
complex-valued harmonic functions, which is as follows: If $f$ is a
complex-valued harmonic function from $\mathbb{D}$ into itself,
then, for $z\in\mathbb{D}$,
\be\label{eq-Co}
\|D_{f}(z)\|\leq\frac{4}{\pi}\frac{1}{1-|z|^{2}}.\ee

The following result establishes a Schwarz-Pick type lemma for the solutions to the $\alpha$-harmonic equation.
\begin{thm}\label{thm-2}
Suppose that $g\in\mathcal{C}(\overline{\mathbb{D}})$, $f\in\mathcal{C}^{2}(\mathbb{D})$ satisfies \eqref{eq1.2} with $\alpha\geq0$ and that $f^{*}\in C(\mathbb{T})$. Then
for $z\in\mathbb{D}$,
$$
\|D_{f}(z)\|\leq(\alpha+1)2^{\alpha+1}\|f^{\ast}\|_{\infty}\frac{1}{1-|z|^2}+(\alpha+\frac{4}{3})2^{\alpha+1}\|g\|_{\infty},
$$
where $\|f^{\ast}\|_{\infty}$ and  $\|g\|_{\infty}$ are as in Theorem \ref{thm-1}.

In particular, if $f$ maps $\ID$ into $\ID$ and $g\equiv0$, then
\be\label{eq-2-1}\|D_{f}(z)\|\leq(\alpha+1)2^{\alpha+1}\frac{1}{1-|z|^2}.\ee
\end{thm}

We remark that the estimate of \eqref{eq-2-1} is sharper than the estimate given by \cite[Theorem 1.1]{LWX}.

\subsection{Compositions of $\alpha$-harmonic functions}
Although the composition $f\circ \phi$, where $f$ is harmonic and $\phi$ is analytic, is known to be harmonic, an $\alpha$-harmonic function ($\alpha\neq0$) precomposed with an analytic function may not be $\alpha$-harmonic. This can be seen by the following example.

\beg\label{ex1}
Let $\alpha\in(-1,0)\cup(0,+\infty)$ and let $k\geq1$ be an integer. Denote by $f$ the function
\beqq f(z)&=&\Big(\sum^{\infty}_{n=0}\frac{\Gamma(n+k)}{\Gamma(n+1)\Gamma(k)}(1-|z|^2)^n\\
&&-(1-|z|^2)^{\alpha+1}\sum^{\infty}_{n=0}\frac{\Gamma(n+\alpha+k+1)}{\Gamma(n+\alpha+2)\Gamma(k)}(1-|z|^2)^n\Big)\overline{z}^k\eeqq
for $z\in\ID\setminus\{0\}$, where $\Gamma(s)=\int^{\infty}_{0}t^{s-1}e^{-t}dt$ $(s>0)$ is the Gamma function. Let $\psi(z)=z^2$. Then
\begin{enumerate}
\item
$f$ is an $\alpha$-harmonic function in $\ID\setminus\{0\}$;
\item
$f\circ \psi$ is not an $\alpha$-harmonic function in $\ID\setminus\{0\}$.
\end{enumerate}

To see this note that, by \cite[Lemma 1.6]{oA}, $f$ is an $\alpha$-harmonic function in $\ID\setminus\{0\}$.
By letting $\eta=\psi(z)=z^2$, we have
$$
\frac{\partial}{\partial \overline{z}}f\circ\psi(z)=\frac{\partial}{\partial \overline{\eta}}f(\eta)\cdot\frac{\partial \overline{\eta}}{\partial\overline{z}}
=2\frac{\Gamma(\alpha+k+1)}{\Gamma(\alpha+1)(k-1)!}(1-|z|^{2})^{\alpha}(1+|z|^{2})^{\alpha}{\overline{z}}^{2k-1},
$$
and then
$$\frac{\partial}{\partial z}\Big((1-|z|^{2})^{-\alpha}\frac{\partial}{\partial \overline{z}}f\circ\psi(z)\Big)\neq0.$$
Hence, it follows that $f\circ \psi$ is not an $\alpha$-harmonic function in $\ID\setminus\{0\}$.
\eeg

%\qed \medskip

%xxxx

Now we are ready to state our result.

\begin{thm}\label{thm-3}
Let $\psi$ be an analytic function in $\ID$. Then for any $\alpha$-harmonic function $f$
with $\alpha\in(-1,0)\cup(0,+\infty)$, $f\circ \psi$ is $\alpha$-harmonic if and only if $\psi(z)=e^{it}z$, $t\in[0,2\pi]$.
\end{thm}

\subsection{Bergman-type spaces}
For $\nu, \mu, t\in\mathbb{R}$,
$$\mathcal{D}_f(\nu, \mu, t)=\int_{\ID}d^{\nu}|f(z)|^{\mu}\|D{f}(z)\|^tdA(z)<\infty$$
is called {\it Dirichlet-type energy integral} of the complex-valued function $f$ (cf. \cite{A, AIM, CRW, CPR, CRV, E, GPP, GP, S}).
In particular, for $\nu>-1$, $0<\mu<\infty$ and $t=0$, we denote by $b_{\nu, \mu}(\ID)$ the {\it Bergman-type space}, consisting of all $f\in\mathcal{C}^0(\ID)$
with the norm
$$\|f\|_{b_{\nu, \mu}}=|f(0)|+(\mathcal{D}_f(\nu, \mu, 0))^{\frac{1}{\mu}}<\infty.$$

We refer to \cite{CRV, DS, GPP, GP, zhu} for basic characterizations of analytic
or harmonic Bergman-type spaces and Dirichlet-type spaces. However, very few related studies can be found from the literature for the
general complex-valued functions.
The following is a characterization of $\alpha$-harmonic functions in Bergman-type spaces.

\begin{thm}\label{thm-4}
Let $f\in\mathcal{C}^{2}(\mathbb{D})$ be an $\alpha$-harmonic function in $\ID$ with $\alpha>-1$, $\real(f\overline{\Delta f})\geq0$ and $\sup_{z\in\overline{\mathbb{D}}}|f(z)|\leq M$, where $M$ is a constant. Then for $p\geq2$, $f\in b_{p-1, p}(\ID)$.
\end{thm}

We will give proofs of Theorem \ref{thm-1} and Theorem \ref{thm-2} in Section \ref{sec-2}.  Theorem \ref{thm-3} and Theorem \ref{thm-4} will be
proved in Section \ref{sec-3}.

\section{Schwarz type and Schwarz-Pick type lemmas for the solutions of the $\alpha$-harmonic equation}\label{sec-2}

The aim of this section is to prove Theorems \ref{thm-1} and \ref{thm-2}. First, we show Theorem \ref{thm-1}. For this, we need some lemmas.

\begin{lem}\label{lem1}
For $\alpha\geq0$, the function $h(s)$ given by \eqref{eqhs} satisfies the estimate
$$h(s)\leq s^{\alpha}\log\frac{1}{1-s}.$$
\end{lem}

\bpf
By \eqref{eqhs}, we obtain that
$$
h(s)=\sum^{\infty}_{n=0}\frac{s^{\alpha+1+n}}{\alpha+1+n}=\sum^{\infty}_{n=1}\frac{s^{\alpha+n}}{\alpha+n}=s^{\alpha}\sum^{\infty}_{n=1}\frac{s^{n}}{\alpha+n}
\leq s^{\alpha}\sum^{\infty}_{n=1}\frac{s^{n}}{n}=s^{\alpha}\log\frac{1}{1-s},$$
and the result follows.
\epf

\begin{lem}\label{lem2}
For $\alpha\geq0$, the function $G_{\alpha}(z,w)$ given by \eqref{eqGr} satisfies the estimate
$$|G_{\alpha}(z,w)|\leq 2^{\alpha}(1-|z|^2)^{\alpha}\log\Big|\frac{1-\overline{z}w}{z-w}\Big|^2.$$
\end{lem}
\bpf
Applying Lemma \ref{lem1}, shows that
\beqq
|G_{\alpha}(z,w)|&\leq& |1-z\overline{w}|^{\alpha}\frac{(1-|z|^2)^{\alpha}(1-|w|^2)^{\alpha}}{|1-z\overline{w}|^{2\alpha}}\log\Big|\frac{1-\overline{z}w}{z-w}\Big|^2\\
&\leq& 2^{\alpha}(1-|z|^2)^{\alpha}\log\Big|\frac{1-\overline{z}w}{z-w}\Big|^2,
\eeqq
and the result follows.
\epf

\subsection*{ Proof of  Theorem \ref{thm-1}}
For a given $g\in\mathcal{C}(\overline{\mathbb{D}})$, by \eqref{eq1.4}, we have
$$|f(z)|\leq|\mathcal{P}_{\alpha}[f^{\ast}](z)|+|\mathcal{G}[g](z)|.$$

Since
\beqq
|\mathcal{P}_{\alpha}[f^{\ast}](z)|&=&\left|\frac{1}{2\pi}\int^{2\pi}_{0}\frac{1-|z|^{2}}{|1-ze^{-i\theta}|^2}\cdot
\frac{(1-|z|^{2})^{\alpha}}{(1-\overline{z}e^{i\theta})^{\alpha}}f^{\ast}(e^{i\theta})d\theta\right|\\
&\leq&(1+|z|)^{\alpha}\frac{1}{2\pi}\int^{2\pi}_{0}\frac{1-|z|^{2}}{|1-ze^{-i\theta}|^2}|f^{\ast}(e^{i\theta})|d\theta\\
&\leq&2^{\alpha}\mathcal{P}[|f^{\ast}|](z)
\eeqq
and $\mathcal{P}[|f^{\ast}|](z)$ is harmonic in $\ID$, by \eqref{eqh}, we know that, for $z\in\ID$
\be\label{eq-P}
|\mathcal{P}_{\alpha}[f^{\ast}](z)|\leq 2^{\alpha}\mathcal{P}[|f^{\ast}|](z)\leq 2^{\alpha}\cdot\frac{4}{\pi}\|f^{\ast}\|_{\infty}\arctan|z|.
\ee

On the other hand, by Lemma \ref{lem2}, we get
$$|\mathcal{G}[g](z)|=\Big|\int_{\ID}G_{\alpha}(z, w)g(w)\,d A(w)\Big|
\leq 2^{\alpha}\|g\|_{\infty}(1-|z|^{2})^{\alpha}J_1,$$
where
$$J_1=\int_{\ID}\log\Big|\frac{1-\overline{z}w}{z-w}\Big|^2\,d A(w).$$

In order to estimate $J_1$, we let
$$w \mapsto \zeta =\phi (w)=\frac{z-w}{1-w\overline{z}}=re^{i\vartheta}
$$
so that $\phi =\phi^{-1}$,
$$w=\frac{z-\zeta}{1-\zeta\overline{z}}, ~~\phi ' (w)=-\frac{1-|z|^2}{(1-w\overline{z})^2},
$$
and thus,
$$d A(w)= |({\phi^{-1}}) ' (\zeta)|^2 d A(\zeta)=\frac{(1-|z|^2)^2}{|1-\zeta\overline{z}|^4}d A(\zeta).
$$
Consequently, by switching to the polar coordinates and using \cite[(2.3)]{LP}, we obtain
\beqq
J_1&=&\int_{\ID}\frac{(1-|z|^2)^2}{|1-\zeta\overline{z}|^4}\log\frac{1}{|\zeta|^2}\,d A(\zeta)\\
&=&\frac{(1-|z|^2)^2}{\pi}\int^{1}_{0}\int^{2\pi}_{0}\frac{r}{|1-\overline{z}re^{i\vartheta}|^4}\log\frac{1}{r^2}\,d\vartheta\, dr\\
&=& 2(1-|z|^2)^2\sum^{\infty}_{n=0}(n+1)^2|z|^{2n}\int^{1}_{0}r^{2n+1}\log\frac{1}{r^2}\, dr\\
&=&1-|z|^2,
\eeqq
which implies that
\be\label{eq-G}|\mathcal{G}[g](z)|\leq 2^{\alpha}\|g\|_{\infty}(1-|z|^{2})^{\alpha+1}.\ee
Hence, it follows from \eqref{eq-P} and \eqref{eq-G} that \eqref{eqs} holds, and the proof of the theorem is complete.

\qed\medskip

In order to prove Theorem \ref{thm-2}, we need some auxiliary lemmas. The following result is from \cite{LWX}.

\begin{Lem}$($\cite[Lemma 2.1]{LWX}$)$\label{lem10}
If $\alpha>-1$ and $f^{*}\in C(\mathbb{T})$, then
$$\frac{\partial}{\partial z}\int^{2\pi}_{0}P_{\alpha}(ze^{-i\theta})f^{\ast}(e^{i\theta})d\theta
=\int^{2\pi}_{0}\frac{\partial}{\partial z}P_{\alpha}(ze^{-i\theta})f^{\ast}(e^{i\theta})d\theta$$
and
$$\frac{\partial}{\partial \overline{z}}\int^{2\pi}_{0}P_{\alpha}(ze^{-i\theta})f^{\ast}(e^{i\theta})d\theta
=\int^{2\pi}_{0}\frac{\partial}{\partial \overline{z}}P_{\alpha}(ze^{-i\theta})f^{\ast}(e^{i\theta})d\theta.$$
\end{Lem}
\medskip

\begin{lem}\label{lem3}
Assume that $f^{*}\in C(\mathbb{T})$ and $\alpha\geq0$. Then
$$\|D_{\mathcal{P}_{\alpha}[f^{\ast}]}(z)\|\leq(\alpha+1)2^{\alpha+1}\|f^{\ast}\|_{\infty}\frac{1}{1-|z|^2}.$$
\end{lem}

\bpf
By elementary calculations, we obtain
$$\frac{\partial}{\partial z}P(ze^{-i\theta})=\frac{e^{-i\theta}}{(1-ze^{-i\theta})^2}$$ and
$$\frac{\partial}{\partial \overline{z}}P(ze^{-i\theta})=\frac{e^{i\theta}}{(1-\overline{z}e^{i\theta})^2}.$$
Then
\beq\label{eq31}\frac{\partial}{\partial{z}}P_{\alpha}(ze^{-i\theta})
&=&\frac{(1-|z|^{2})^{\alpha}\left[e^{-i\theta}(1-|z|^{2})-(\alpha+1)\overline{z}(1-ze^{-i\theta})\right]}
{(1-ze^{-i\theta})^{2}(1-\overline{z}e^{i\theta})^{\alpha+1}}\nonumber\\
&=&\frac{(1-|z|^{2})^{\alpha}}{(1-\overline{z}e^{i\theta})^{\alpha+1}}\big[1-|z|^{2}-(\alpha+1)\overline{z}e^{i\theta}(1-ze^{-i\theta})\big]
\frac{\partial}{\partial z}P(ze^{-i\theta})
\eeq
and
\beq\label{eq32}
\frac{\partial}{\partial{\overline{z}}}P_{\alpha}(ze^{-i\theta})=\frac{(\alpha+1)(1-|z|^{2})^{\alpha}e^{i\theta}}{(1-\overline{z}e^{i\theta})^{\alpha+2}}
=\frac{(\alpha+1)(1-|z|^{2})^{\alpha}}{(1-\overline{z}e^{i\theta})^{\alpha}}\frac{\partial}{\partial \overline{z}}P(ze^{-i\theta}).
\eeq

Since
\beqq
\big|1-|z|^{2}-(\alpha+1)\overline{z}e^{i\theta}(1-ze^{-i\theta})\big|
&=&\big|(1-ze^{-i\theta})(-\alpha \overline{z}e^{i\theta})+(1-\overline{z}e^{i\theta})\big|\\
&\leq&(\alpha+1)|1-\overline{z}e^{i\theta}|,
\eeqq
we see from \eqref{eq31} that
$$\left|\frac{\partial}{\partial{z}}P_{\alpha}(ze^{-i\theta})\right|
\leq(\alpha+1)2^{\alpha}\left|\frac{\partial}{\partial{z}}P(ze^{-i\theta})\right|.$$
Hence, by combining the above with \eqref{eq32}, we conclude
$$\left|\frac{\partial}{\partial{z}}P_{\alpha}(ze^{-i\theta})\right|+\left|\frac{\partial}{\partial{\overline{z}}}P_{\alpha}(ze^{-i\theta})\right|
\leq(\alpha+1)2^{\alpha}\left(\left|\frac{\partial}{\partial{z}}P(ze^{-i\theta})\right|+\left|\frac{\partial}{\partial{\overline{z}}}P(ze^{-i\theta})\right|\right).$$
Consequently, from Lemma \Ref{lem10} and the identities
$$\frac{1}{2\pi}\int^{2\pi}_{0}\left|\frac{\partial}{\partial{z}}P(ze^{-i\theta})\right|d\theta
=\frac{1}{2\pi}\int^{2\pi}_{0}\left|\frac{\partial}{\partial{\overline{z}}}P(ze^{-i\theta})\right|d\theta=\frac{1}{1-|z|^{2}},$$
it follows that
$$\|D_{\mathcal{P}_{\alpha}[f^{\ast}]}(z)\|\leq(\alpha+1)2^{\alpha+1}\|f^{\ast}\|_{\infty}\frac{1}{1-|z|^2},$$
which is what is needed.
\epf

Identify the complex plane $\mathbb{C}$ with $\mathbb{R}^2$, and denote by $L_{loc}^1(\ID)$ the space of locally integrable functions in $\ID$.
Functions $\psi\in L_{loc}^1(\ID)$ with distributional partial derivatives in $\ID$, we have the action
$$
\langle\psi_z, \varphi\rangle=-\int_{\ID}\psi\varphi_z \,dA,\;\;\varphi\in C_0^{\infty}(\ID),
$$
for the distributional partial derivative $\psi_z$ have, and similarly for the distribution $\psi_{\overline{z}}$ (cf. \cite{B, GG, LH}).

%hav the products
%$$\langle\psi, \varphi\rangle=\int_{\ID}\psi\varphi \,dA,\;\;\varphi\in C_0^{\infty}(\ID).
%$$

\begin{lem}\label{lem4}
For $\alpha\geq0$ and $g\in\mathcal{C}(\overline{\mathbb{D}})$, the function $G_{\alpha}(z,w)$, given by \eqref{eqGr}, satisfies the following inequality:
\bee
\item $\ds \int_{\ID}\left|G_{\alpha}(z,w)g(w)\right|d A(w)\leq 2^{\alpha}\|g\|_{\infty}$ and $\int_{\ID}\left|G_{\alpha}(z, w)\right|d A(z)<\infty$;
\item For fixed $w\in\ID$,
$$
\frac{\partial G_{\alpha}(z,w)}{\partial z}=\alpha \overline{w}(1-z\overline{w})^{\alpha-1}h\circ g(z, w)
-\frac{(1-|z|^{2})^{\alpha}(1-|w|^{2})^{\alpha+1}}{(1-\overline{z}w)^{\alpha}(1-z\overline{w})}\cdot\frac{1}{z-w}
$$
and
$$\frac{\partial G_{\alpha}(z,w)}{\partial \overline{z}}=\frac{(1-|z|^{2})^{\alpha}(1-|w|^{2})^{\alpha+1}}{(1-\overline{z}w)^{\alpha+1}(\overline{z}-\overline{w})}$$
in the sense of distributions in $\ID;$
\item
\bee
\item $\ds \frac{\partial \mathcal{G}[g](z)}{\partial z}=\int_{\ID}\frac{\partial G_{\alpha}(z, w)}{\partial z}g(w)\,d A(w)$,
and
\item $\ds \left| \frac{\partial \mathcal{G}[g](z)}{\partial z}\right|\leq \int_{\ID}\left|\frac{\partial G_{\alpha}(z, w)}{\partial z}g(w)\right|d A(w)
\leq (\alpha+\frac{2}{3})2^{\alpha+1}\|g\|_{\infty}$,
\eee
in the sense of distributions in $\ID;$
\item
\bee
\item
$\ds \frac{\partial \mathcal{G}[g](z)}{\partial \overline{z}}=\int_{\ID}\frac{\partial G_{\alpha}(z, w)}{\partial \overline{z}}g(w)\,d A(w)$,
and
\item
$\ds \left| \frac{\partial \mathcal{G}[g](z)}{\partial \overline{z}}\right|\leq\int_{\ID}\left|\frac{\partial G_{\alpha}(z, w)}{\partial \overline{z}}g(w)\right|d A(w)
\leq \frac{2^{\alpha+2}}{3}\|g\|_{\infty}$,
\eee
in the sense of distributions in $\ID$.
\eee
\end{lem}

\bpf
Obviously, it follows from \eqref{eq-G} and the proof of \cite[Proposition 4]{B} that $(1)$ holds. Therefore, $G_{\alpha}(z, w)\in L^1(\ID)$ so that its derivative has the action
$$\left\langle \frac{\partial G_{\alpha}(z, w)}{\partial z}, \varphi(z)\right\rangle=-\int_{\ID} G_{\alpha}(z, w)\varphi_{z}(z)\, d A(z), \;\;\varphi\in C_0^{\infty}(\ID).
$$
By Lebesgue's dominated convergence theorem we get
$$\int_{\ID} G_{\alpha}(z, w)\varphi_{z}(z)\, d A(z)
=\lim_{\varepsilon\to0}\int_{\ID\setminus \mathbb{D}(w, \varepsilon)} G_{\alpha}(z, w)\varphi_{z}(z) \, d A(z).
$$
For $\varepsilon>0$, let $D_{\varepsilon}=D(w, \varepsilon)$. Partial integration gives
\beqq
\int_{\ID\setminus D_{\varepsilon}} G_{\alpha}(z, w)\varphi_{z}(z) \,d A(z)&=&-\int_{\partial D_{\varepsilon}} G_{\alpha}(z, w)\varphi(z)v(z)\, d s(z)\\
&&-\int_{\ID\setminus D_{\varepsilon}}\frac{\partial G_{\alpha}(z, w)}{\partial z}\varphi(z)\,d A(z),
\eeqq
where $v$ is the unit exterior normal of $\ID\setminus D_{\varepsilon}$, that is, the inward directed unit normal of $D_{\varepsilon}$ and $ds$ denotes the normalized
arc length measure.

It follows from Lemma \ref{lem2} that
$$|G_{\alpha}(z,w)|\leq 2^{\alpha}(1-|z|^2)^{\alpha}\log\Big|\frac{1-\overline{z}w}{z-w}\Big|^2.$$
Then
\beqq
\left|\int_{\partial D_{\varepsilon}} G_{\alpha}(z,w)\varphi(z)v(z)\, d s(z)\right|
&\leq& C_{\varphi}\int_{\partial D_{\varepsilon}} \left|G_{\alpha}(z,w)\right|d s(z)\\
&\leq& 2^{\alpha}C_{\varphi}\varepsilon\log \frac{4}{{\varepsilon}^2}\rightarrow0,
\eeqq
as $\varepsilon\rightarrow0$, where $C_{\varphi}$ is a constant depending only on $\sup\varphi$, and the supremum is taken over all functions $\varphi$.
A straightforward computation shows that for $z\neq w$ we have
$$
\frac{\partial G_{\alpha}(z,w)}{\partial z}=H(z, w)=\alpha \overline{w}(1-z\overline{w})^{\alpha-1}h\circ g(z, w)
+\frac{(1-|z|^{2})^{\alpha}(1-|w|^{2})^{\alpha}}{(1-\overline{z}w)^{\alpha}(1-z\overline{w})}\cdot\frac{1}{z-w}.
$$
It follows from Lemma \ref{lem1} that
\beq\label{eq41}
\left|\frac{\partial G_{\alpha}(z, w)}{\partial z}\right|&\leq& \alpha\cdot2^{\alpha+1}(1-|z|^{2})^{\alpha-1}\log\Big|\frac{1-\overline{z}w}{z-w}\Big|^2\\
&&+2^{\alpha}(1-|z|^{2})^{\alpha}\frac{1-|w|^{2}}{|1-z\overline{w}|\cdot|z-w|}.\nonumber
\eeq
Therefore
\beqq
\left\langle \frac{\partial G_{\alpha}(z, w)}{\partial z}, \varphi(z)\right\rangle
&=&\lim_{\varepsilon\to0}\int_{\ID\setminus D_{\varepsilon}}\frac{\partial G_{\alpha}(z, w)}{\partial z}\varphi(z)\, d A(z)\\
&=&\int_{\ID}H(z,w)\varphi(z)\, d A(z)
=\langle H(z,w), \varphi(z)\rangle.
\eeqq
By a similar reasoning as above, one obtains that
$$\frac{\partial G_{\alpha}(z,w)}{\partial \overline{z}}=\frac{(1-|z|^{2})^{\alpha}(1-|w|^{2})^{\alpha+1}}{(1-\overline{z}w)^{\alpha+1}(\overline{z}-\overline{w})}$$
in the sense of distributions in $\ID$. Hence, the assertion (2) of the lemma holds.

To prove the assertion (3),  it follows from \eqref{eq41} that
\beqq
\int_{\ID}\left|\frac{\partial G_{\alpha}(z, w)}{\partial z}\right|d A(w)
&\leq&\alpha2^{\alpha+1}(1-|z|^{2})^{\alpha-1}\int_{\ID}\log\left|\frac{1-\overline{w}z}{z-w}\right|^2\,d A(w)\\
&&+2^{\alpha}(1-|z|^{2})^{\alpha}\int_{\ID}\frac{1-|w|^{2}}{|1-z\overline{w}|\cdot|z-w|}\,d A(w).
\eeqq
Moreover, as before, by using the transformation
$$w \mapsto \zeta =\phi (w)=\frac{z-w}{1-w\overline{z}}=re^{i\vartheta}$$
we obtain after computation that
\beqq
\int_{\ID}\frac{1-|w|^2}{|z-w|\cdot|1-\overline{w}z|}\, d A(w)
&=&\int_{\ID}\frac{(1-|z|^2)(1-|\zeta|^2)}{|\zeta|\cdot|1-\zeta\overline{z}|^4}\, d A(\zeta)\\
&=&\frac{(1-|z|^2)}{\pi}\int^{1}_{0}\int^{2\pi}_{0}\frac{1-r^2}{|1-\overline{z}re^{i\vartheta}|^4}\,d\vartheta\, dr\\
&=&2(1-|z|^2)\sum^{\infty}_{n=0}(n+1)^2|z|^{2n}\int^{1}_{0}r^{2n}(1-r^2)\, dr \\
&=&4(1-|z|^2)\sum^{\infty}_{n=0}\frac{(n+1)^2}{(2n+1)(2n+3)} |z|^{2n}  \\
&\leq&\frac{4(1-|z|^2)}{3}\sum^{\infty}_{n=0} |z|^{2n} = \frac{4}{3}. %\quad (\mbox{since $3(n+1)^2\leq (2n+1)(2n+3)$})\\
\eeqq
Thus, we conclude that
\beqq
\int_{\ID}\left|\frac{\partial G_{\alpha}(z, w)}{\partial z}\right|d A(w)
&\leq&\alpha2^{\alpha+1}(1-|z|^{2})^{\alpha-1}J_1+2^{\alpha}\frac{4}{3}(1-|z|^{2})^{\alpha}\\
&=&\alpha2^{\alpha+1}(1-|z|^{2})^{\alpha}+2^{\alpha}\frac{4}{3}(1-|z|^{2})^{\alpha}\\
&\leq&(\alpha+\frac{2}{3})2^{\alpha+1}.
\eeqq
Hence
\beq\label{eq42}
\int_{\ID}\left|\frac{\partial G_{\alpha}(z, w)}{\partial z}g(w)\right|d A(w)
\leq (\alpha+\frac{2}{3})2^{\alpha+1}\|g\|_{\infty}.
\eeq

By the assertion (1) and \eqref{eq42}, we see that
\beqq
\left\langle \frac{\partial }{\partial z}\int_{\ID}G_{\alpha}(z,w)g(w)d A(w), \varphi(z)\right\rangle
&=&-\int_{\ID} \left(\int_{\ID}G_{\alpha}(z,w)g(w)d A(w)\right)\varphi_{z}(z)\, d A(z)\\
&=&-\int_{\ID} \left(\int_{\ID}G_{\alpha}(z,w)g(w)\varphi_{z}(z)d A(z)\right)\, d A(w)\\
&=&\int_{\ID} \left\langle \frac{\partial }{\partial z}G_{\alpha}(z,w)g(w), \varphi(z)\right\rangle\, d A(w)\\
&=&\int_{\ID} \left(\int_{\ID}\frac{\partial }{\partial z}G_{\alpha}(z,w)g(w)\varphi(z)d A(z)\right)\, d A(w)\\
&=&\int_{\ID} \left(\int_{\ID}\frac{\partial }{\partial z}G_{\alpha}(z,w)g(w)d A(w)\right)\varphi(z)\, d A(z)\\
&=&\left\langle \int_{\ID}\frac{\partial }{\partial z}G_{\alpha}(z,w)g(w)d A(w), \varphi(z)\right\rangle
\eeqq
Hence, we conclude that
$$\ds \frac{\partial \mathcal{G}[g](z)}{\partial z}=\int_{\ID}\frac{\partial G_{\alpha}(z, w)}{\partial z}g(w)\,d A(w),$$
and then the assertion (3) holds.

For the assertion (4), it follows from the assertion (2) that
$$\left|\frac{\partial G_{\alpha}(z,w)}{\partial \overline{z}}\right|
\leq 2^{\alpha}(1-|z|^{2})^{\alpha}\frac{1-|w|^{2}}{|1-\overline{z}w||z-w|},$$
and then
\beqq
\int_{\ID}\left|\frac{\partial G_{\alpha}(z, w)}{\partial \overline{z}}g(w)\right|d A(w)
&\leq&2^{\alpha}(1-|z|^{2})^{\alpha}\|g\|_{\infty}\int_{\ID}\frac{1-|w|^2}{|z-w|\cdot|1-\overline{w}z|}\, d A(w)\\
&\leq&\frac{2^{\alpha+2}}{3}\|g\|_{\infty}.
\eeqq
The proof of the remaining assertion (4) is similar to the proof of the assertion (3).
Hence the proof of the lemma is complete.
\epf

\subsection*{ Proof of  Theorem \ref{thm-2}}

The result follows from Lemma \ref{lem3} and Lemma \ref{lem4} together with \eqref{qh-1}.
\qed
\medskip

\section{Some properties of $\alpha$-harmonic functions}\label{sec-3}
In this section, we will prove Theorem \ref{thm-3} and Theorem \ref{thm-4}.

% First, we give the properties of Example \ref{ex1}.

In \cite{oA}, the authors obtained the following homogeneous expansion of $\alpha$-harmonic functions
(see \cite[Theorem 1.2]{oA}):

A function $f$ in $\mathbb{D}$ is $\alpha$-harmonic if and only if it has the following convergent power series expansion:
\beq\label{eq1.31}\;\;\;
f(z)=\sum^{\infty}_{k=0}c_{k}z^{k}+\sum^{\infty}_{k=1}c_{-k}P_{\alpha,k}(|z|^2)\overline{z}^{k},
\eeq
where
\beq\label{eq1.32}
P_{\alpha,k}(x)=\int^{1}_{0}t^{k-1}(1-tx)^{\alpha}dt,
\eeq
 $-1<x<1$, and $\{c_{k}\}^{\infty}_{k=-\infty}$ denotes a sequence of complex numbers with
$$
\lim_{|k|\to\infty}\sup|c_{k}|^{\frac{1}{|k|}}\leq1.
$$

Now, we are ready to prove Theorem \ref{thm-3}.

\subsection*{ Proof of  Theorem \ref{thm-3}}
The sufficiency is obvious, and now we prove the necessity.
Suppose that $f\circ\psi$ is $\alpha$-harmonic. By \eqref{eq1.31}, we get
$$f\circ\psi(z)=\sum^{\infty}_{k=0}c_{k}\psi(z)^{k}+\sum^{\infty}_{k=1}c_{-k}P_{\alpha,k}(|\psi(z)|^2)\overline{\psi(z)}^{k}.$$
Take the derivative of  \eqref{eq1.32}, we obtain
$$P_{\alpha,k}'(x)=-\int^{1}_{0}t^{k}\alpha(1-tx)^{\alpha-1}dt,$$
and an integration by parts gives
\beq\label{eq1.33}
xP_{\alpha,k}'(x)+k P_{\alpha,k}(x)=(1-x)^{\alpha}.
\eeq
Let $\tau(z)=P_{\alpha,k}(|\psi(z)|^2)\overline{\psi(z)}^{k}$. Then, by \eqref{eq1.33}, we obtain
\beqq
\tau_{\overline{z}}(z)&=&\overline{\psi'(z)}\cdot\overline{\psi(z)}^{k-1}\Big(|\psi(z)|^2P_{\alpha,k}'(|\psi(z)|^2)+kP_{\alpha,k}(|\psi(z)|^2)\Big)\\
&=&\overline{\psi'(z)}\cdot\overline{\psi(z)}^{k-1}(1-|\psi(z)|^2)^{\alpha}.
\eeqq
This together with the fact that $f\circ\psi$ is $\alpha$-harmonic, shows that
$$\frac{\partial}{\partial z}\Big[(1-|z|^2)^{-\alpha}(1-|\psi(z)|^2)^{\alpha}\Big]=0,$$
and then
\beq\label{eq1.34}
z\psi'(z)\overline{\psi(z)}=\frac{1-|\psi(z)|^2}{1-|z|^2}|z|^2
\eeq
is real valued. This implies that $\psi(z)=a_kz^k$, $k=1, 2,\ldots$, where $a_k$ are constants.
Therefore, it follows from \eqref{eq1.34} that, for $z\in\ID\backslash\{0\}$,
$$|a_k|^2=\frac{1}{k|z|^{2k-2}-(k-1)|z|^{2k}}.$$
Hence, for $k\geq2$, $a_k$ are not constants. This contradicts with the assumption that $a_k$ are constants.
Therefore, $\psi(z)=a_1z$. By \eqref{eq1.34}, we have that
$$|a_1|^2|z|^2=\frac{1-|a_1|^2|z|^2}{1-|z|^2}|z|^2.$$
Therefore, we conclude that $|a_1|^2=1$, which completes the proof.
\qed
\medskip

\subsection*{ Proof of  Theorem \ref{thm-4}}
By \cite[Theorem 3]{CRV}, we only need to prove
\beq\label{eq1.41}
\int_{\ID}(1-|z|^2)^{p+1}\Delta(|f(z)|^p)\, d A(z)<\infty.
\eeq
From Theorem \ref{thm-2} we know that
\beq\label{eq1.42}
\|D_{f}(z)\|\leq \frac{C_1}{1-|z|^2},
\eeq
where $C_1=(\alpha+1)2^{\alpha+1}$.
Since $f$ is $\alpha$-harmonic, the function $(1-|z|^2)^{-\alpha}f_{\overline{z}}$ is antianalytic.
Therefore it has a power series expansion of the form
$$(1-|z|^2)^{-\alpha}f_{\overline{z}}(z)=\sum^{\infty}_{k=0}a_{k}\overline{z}^{k}, \;\;z\in\ID.$$
By calculations, we have
$$\Delta f(z)=4f_{\overline{z}z}(z)=\frac{-4\alpha \overline{z}}{1-|z|^2}f_{\overline{z}}(z).$$
Hence, it follows from \eqref{eq1.42} that
\beq\label{eq1.43}
|\Delta f(z)|\leq\frac{4\alpha}{1-|z|^2}\|D_{f}(z)\|\leq \frac{C_2}{(1-|z|^2)^2},
\eeq
where $C_2=\alpha(\alpha+1)2^{\alpha+3}$. Furthermore,
\beqq
|f(z)|&\leq&|f(0)|+\left|\int_{[0, z]}d f(\xi)\right|\\
&\leq&|f(0)|+\int_{[0, z]}\|D_{f}(\xi)\|\;|d\xi|\\
&\leq&|f(0)|+\frac{C_1}{1-|z|^2},
\eeqq
where $[0, z]$ denotes the line segment from $0$ to $z$.
Then, for $z\in\ID$,
\beq\label{eq1.44}
|f(z)|^{p-1}\leq 2^{p-2}\left(|f(0)|^{p-1}+\frac{C_1^{p-1}}{(1-|z|^2)^{p-1}}\right),
\eeq
and
\beq\label{eq1.45}
|f(z)|^{p-2}\leq 2^{p-2}\left(|f(0)|^{p-2}+\frac{C_1^{p-2}}{(1-|z|^2)^{p-2}}\right).
\eeq

We divide the remaining part of the proof into two cases, namely, $p\in[4,\infty)$ and $p\in[2, 4)$.
For the case $p\in[4,\infty)$, a straightforward computation gives
\beqq
\Delta(|f|^p)&=&p(p-2)|f|^{p-4}|f\overline{f_z}+f_{\overline{z}}\overline{f}|^2+2p|f|^{p-2}(|f_z|^2+|f_{\overline{z}}|^2)+p|f|^{p-2}\real(\overline{f}\Delta f)\\
&\leq& p^2|f|^{p-2}\|D_{f}\|^2+p|f|^{p-1}|\Delta f|.
\eeqq
Hence, by \eqref{eq1.42}, \eqref{eq1.43}, \eqref{eq1.44} and \eqref{eq1.45}, we conclude that
\beqq\label{eq1.46}
&&(1-|z|^2)^{p+1}\Delta(|f(z)|^p)\nonumber \\
&\leq&p^2(1-|z|^2)^{p+1}|f(z)|^{p-2}\|D_{f}(z)\|^2+p(1-|z|^2)^{p+1}|f(z)|^{p-1}|\Delta f(z)|\nonumber \\
&\leq&p^2 2^{p-2}(1-|z|^2)^{p+1}\left(|f(0)|^{p-2}+\frac{C_1^{p-2}}{(1-|z|^2)^{p-2}}\right)\frac{C_1^2}{(1-|z|^2)^2}\nonumber \\
&&+p2^{p-2}(1-|z|^2)^{p+1}\left(|f(0)|^{p-1}+\frac{C_1^{p-1}}{(1-|z|^2)^{p-1}}\right)\frac{C_2}{(1-|z|^2)^2}\nonumber \\
&\leq&p^2 2^{p-2}(|f(0)|^{p-2}+C_1^{p-2})C_1^2(1-|z|^2)+p2^{p-2}(|f(0)|^{p-1}+C_1^{p-1})C_2\nonumber \\
&<&\infty.
\eeqq

In the case $p\in[2, 4)$, let $F_n^p=(|f|^2+\frac{1}{n})^{\frac{p}{2}}$ for $n\in\left\{1, 2, \ldots\right\}$.
Then
\beqq
\Delta(F_n^p)&=&p(p-2)\Big(|f|^2+\frac{1}{n}\Big)^{\frac{p}{2}-2}|f\overline{f_z}+f_{\overline{z}}\overline{f}|^2\\
&&+2p\Big(|f|^2+\frac{1}{n}\Big)^{\frac{p}{2}-1}(|f_z|^2+|f_{\overline{z}}|^2)+p\Big(|f|^2+\frac{1}{n}\Big)^{\frac{p}{2}-1}\real(\overline{f}\Delta f).
\eeqq
Let
\beqq
F&=&p(p-2)|f|^{p-2}\|D_{f}(z)\|^2+2p\big(|f|^2+1\big)^{\frac{p}{2}-1}(|f_z|^2+|f_{\overline{z}}|^2)\\
&&+p\big(|f|^2+1\big)^{\frac{p}{2}-1}\real(\overline{f}\Delta f).
\eeqq
For $r\in(0, 1)$, $\Delta(F_n^p)$ and $F$ are integrable in $\ID_r$, and $\Delta(F_n^p)\leq F$.
By Lebesgue's dominated convergence theorem, we have
\beqq
\lim_{n\rightarrow\infty}\int_{\ID_r}(1-|z|^2)^{p+1}\Delta\big(F_n^p(z)\big)\, d A(z)
&=&\int_{\ID_r}(1-|z|^2)^{p+1}\lim_{n\rightarrow\infty}\Big(\Delta\big(F_n^p(z)\big)\Big)\, d A(z)\\
&\leq&\int_{\ID_r}(1-|z|^2)^{p+1}\Big[p^2|f(z)|^{p-2}\|D_{f}(z)\|^2\\
&&+p|f(z)|^{p-1}|\Delta f(z)|\Big]\, d A(z)\\
&\leq&\int_{\ID_r}\Big[p^2 2^{p-2}(|f(0)|^{p-2}+C_1^{p-2})C_1^2(1-|z|^2)\\
&&+p2^{p-2}(|f(0)|^{p-1}+C_1^{p-1})C_2\Big]\, d A(z)\\
&<&\infty.
\eeqq
Therefore, \eqref{eq1.41} follows from the above two estimates and so the proof of the theorem is complete.
\qed
\medskip

\vskip.10in

\begin{center}\textbf{\sc Acknowledgments}
\end{center}

\vskip.05in

The present investigation was supported by the \textit{Natural
Science Foundation of Hunan Province} under Grant no. 2016JJ2036.

%\medskip

\end{document}